\documentclass[11pt,reqno]{amsart}
\usepackage{graphicx}
\usepackage{verbatim}
\usepackage{textcomp}
\usepackage{amssymb}
\usepackage{cite}
\usepackage{amsmath}
\usepackage{latexsym}
\usepackage{amscd}
\usepackage{amsthm}
\usepackage{mathrsfs}
\usepackage{xypic}
\usepackage{bm}
\usepackage{url}
\usepackage{hyperref}

\newtheorem{thm}{Theorem}[section]

\newtheorem{lem}[thm]{Lemma}

\theoremstyle{definition}

\theoremstyle{remark}
\newtheorem{rem}{Remark}[section]
\numberwithin{equation}{section}
\begin{document}
\title[Hamilton-Souplet-Zhang's gradient estimates]
{Hamilton-Souplet-Zhang's gradient estimates for two types of
nonlinear parabolic equations under the Ricci flow }
\author{Guangyue Huang}
\address{College of Mathematics and Information Science, Henan Normal
University, Xinxiang, Henan 453007, People's Republic of China}
\email{hgy@henannu.edu.cn }
\author{Bingqing Ma}
\address{College of Mathematics and Information Science,
Henan Normal University, Xinxiang, Henan 453007, People's Republic
of China} \email{bqma@henannu.edu.cn} 

\maketitle

\begin{abstract}
In this paper, we consider gradient estimates for two types of
nonlinear parabolic equations under the Ricci flow: one is the
equation
$$u_t=\Delta u+au\log u+bu$$ with $a,b$ two real constants, the
other is $$u_t=\Delta u+\lambda u^{\alpha}$$ with $\lambda,\alpha$
two real constants. By a suitable scaling for the above two
equations, we obtain Hamilton-Souplet-Zhang type gradient estimates.
\end{abstract}

\section{Introduction}

Since the nonlinear parabolic equation
\begin{equation}\label{0Int1}
u_t=\Delta u+au\log u+bu
\end{equation} on a given Riemannian manifold is related to
gradient Ricci solitons which are self-similar solutions to the
Ricci flow, many attentions are paid to the study on gradient
estimates for the equation \eqref{0Int1}, for example, see
\cite{Ma06,Yangyuyan07,Liyau86,huangli13,lixu11,huangma10,huangma15}.
Here $a,b$ in \eqref{0Int1} are two real constants. Clearly, a heat
equation
\begin{equation}\label{0Int2}
u_t=\Delta u
\end{equation}
is a special case of \eqref{0Int1} when $a=b=0$. Hence many known
results on heat equations are generalized to the nonlinear parabolic
equation \eqref{0Int1}. For gradient estimates of solutions to
\eqref{0Int1} of Li-Yau type, Davies type, Hamilton type and Li-Xu
type on a given Riemannian manifold, we refer to
\cite{huangli13,Yangyuyan07,huangma10,Chen09} and the references
therein. In a recent paper \cite{DungKhanh15},  Dung and Khanh
obtained Hamilton-Souplet-Zhang type gradient estimates on a given
Riemannian manifold for \eqref{0Int1}. On a family of Riemannian
metrics $g(t)$ evolving by the Ricci flow
\begin{equation}\label{1Int3}
\frac{\partial}{\partial t}g_{ij}=-2R_{ij},
\end{equation} Hsu in \cite{Hsu11} obtained Li-Yau type gradient estimates of
\eqref{0Int1}.

In  \cite{SoupletZhang2006}, generalizing Hamilton's estimate in
\cite{hamilton93}, Souplet and Zhang proved

\vspace{2mm}

\noindent {\bf Theorem A\cite{SoupletZhang2006}.} {\it Let $(M^n,
g)$ be an $n$-dimensional Riemannian manifold with ${\rm
Ric}(M^n)\geq-k$, where $k$ is a non-negative constant. Suppose that
$u$ is a positive solution to the equation \eqref{0Int2} in
$Q_{R,T}=\{(x,t)\,|\, x\in M,\, d(x,x_0,t)<R,\, t\in [0,T]\}$ with
$u\leq A$. Then in $Q_{R,T}$,
\begin{equation}\label{Zhang}
\frac{|\nabla u|}{u}\leq C\left(
\sqrt{k}+\frac{1}{R}+\frac{1}{\sqrt{T}}\right)\left(1+\log\frac{A}{u}\right),
\end{equation}
where the constant $C$ depends only on the dimension $n$.

}

The key to prove Theorem A of Souplet and Zhang is the scaling
$u\rightarrow \tilde{u}=u/A$. After this scaling, \eqref{0Int2}
becomes the following heat equation with respect to $\tilde{u}$:
\begin{equation}\label{0Int3}
{\tilde{u}}_t=\Delta \tilde{u}
\end{equation}
since the heat equation is linear. Under this case, we obtain that
$0<\tilde{u}\leq 1$. Inspired by this method, in this paper, we also
adopt the similar scaling method by $u\rightarrow \tilde{u}=u/A$ to
study the nonlinear parabolic equation \eqref{0Int1}. By the
scaling, we can derive from \eqref{0Int1} the following analogous
equation:
\begin{equation}\label{0Int4}
\tilde{u}_t=\Delta \tilde{u}+a\tilde{u}\log
\tilde{u}+\tilde{b}\tilde{u},
\end{equation} where the constant $\tilde{b}$ satisfies $\tilde{b}=b+a\log A$.
That is, we only need to study the nonlinear equation \eqref{0Int4}
with $0<\tilde{u}\leq 1$.

Our first result is the following Hamilton-Souplet-Zhang type
gradient estimates of the nonlinear equation \eqref{0Int1} under the
Ricci flow:

\begin{thm}\label{thmInt1}
Let $M$ be a complete Riemannian manifold with a family of
Riemannian metric $g(t)$ evolving by the Ricci flow \eqref{1Int3}.
Suppose that $u$ is a positive solution to \eqref{0Int1} in
\begin{equation}\label{1thm1}
B_{R,T}=\{(x,t)\,|\, x\in M,\, d(x,x_0,t)<R,\, t\in [0,T]\}
\end{equation} with $|{\rm Ric}|\leq k$ for some positive constant
$k$ and $u\leq A$. Then there exists a constant $C$ depending only
on the dimension of $M$ such that
\begin{equation}\label{1thmFormula1}
\frac{|\nabla u|}{u} \leq
C\Bigg(\frac{1}{R}+\frac{1}{\sqrt{T}}+\sqrt{k}+\sqrt{M_{a,b}}\Bigg)\Bigg(1+\log\frac{A}{u}\Bigg)
\end{equation}
for all $(x,t) \in B_{\frac{R}{2},T}$ with $t\neq 0$, where
$M_{a,b}=\max\{0,a(1+\log A)+b\}$.

\end{thm}

The study to Li-Yau type estimates of the following nonlinear
parabolic equation
\begin{equation}\label{1Int2}
u_t=\Delta u+\lambda u^{\alpha},
\end{equation} where $\lambda,\alpha$ are two real constants, can be traced back to Li
\cite{Li91}. Later, for $0<\alpha<1$, Zhu in \cite{Zhu11Nonlinear}
obtained Hamilton-Souplet-Zhang type gradient estimates of
\eqref{1Int2}  on a given Riemannian manifold.
On gradient estimates of the elliptic case of
\eqref{1Int2}, see \cite{Zhang2011,Yang2010}. A natural subject is
to study Hamilton-Souplet-Zhang type gradient estimates of the
nonlinear equation \eqref{1Int2} under the Ricci flow. Our second
result is the following:

\begin{thm}\label{thmInt2}
Let $M$ be a complete Riemannian manifold with a family of
Riemannian metric $g(t)$ evolving by the Ricci flow \eqref{1Int3}.
Suppose that $u$ is a positive solution to \eqref{1Int2} in
\begin{equation}\label{2thm1}
B_{R,T}=\{(x,t)\,|\, x\in M,\, d(x,x_0,t)<R,\, t\in [0,T]\}
\end{equation} with $|{\rm Ric}|\leq k$ for some positive constant
$k$ and $u\leq B$. Then there exists a constant $C$ depending only
on the dimension of $M$ such that

1) if $\alpha\geq1$, then
\begin{equation}\label{2thmFormula1}
\frac{|\nabla u|}{u} \leq
C\Bigg(\frac{1}{R}+\frac{1}{\sqrt{T}}+\sqrt{k}+\sqrt{M_{\lambda}
\alpha}\Bigg)\Bigg(1+\log\frac{B}{u}\Bigg),
\end{equation}
for all $(x,t) \in B_{\frac{R}{2},\frac{T}{2}}$ with $t\neq 0$,
where $M_{\lambda}=\max\{0,\lambda B^{\alpha-1}\}$;

2) if $\alpha\leq0$, then
\begin{equation}\label{2thmFormula2}
\frac{|\nabla u|}{u} \leq
C\Bigg(\frac{1}{R}+\frac{1}{\sqrt{T}}+\sqrt{k}
+\sqrt{M_{\lambda}(-\alpha+1)u_{\min}^{\alpha-1}}\Bigg)\Bigg(1+\log\frac{B}{u}\Bigg),
\end{equation}
for all $(x,t) \in B_{\frac{R}{2},T}$ with $t\neq 0$, where
$M_{\lambda}=\max\{0,-\lambda \}$;

3) if $\alpha\in(0,1)$, then
\begin{equation}\label{2thmFormula3}
\frac{|\nabla u|}{u} \leq
C\Bigg(\frac{1}{R}+\frac{1}{\sqrt{T}}+\sqrt{k}
+\sqrt{|\lambda|u_{\min}^{\alpha-1}}\Bigg)\Bigg(1+\log\frac{B}{u}\Bigg)
\end{equation}
for all $(x,t) \in B_{\frac{R}{2},T}$ with $t\neq 0$.

\end{thm}

\begin{rem}
Taking $a=b=0$ in \eqref{1thmFormula1}, we obtain the estimate (2.3)
of Theorem 2.2 in \cite{Cao2010} with respect to the heat equation
under the Ricci flow. Hence, our estimates in Theorem \ref{thmInt1}
extend Bailesteanu, Cao and Pulemotov's Theorem 2.2.
\end{rem}

\begin{rem}
There are many studies on gradient estimates of the heat equation
\eqref{0Int2} under geometric flows, we refer to \cite{Liu09,Sun11}
and among others.
\end{rem}


\section{Proof of Theorem \ref{thmInt1}}

In order to prove our Theorem \ref{thmInt1}, we first give a lemma
which will play an important role in the proof.

\begin{lem}\label{2Lem1}
Let $M$ be a complete Riemannian manifold with a family of
Riemannian metric $g(t)$ evolving by the Ricci flow \eqref{1Int3}.
Let $u$ be a positive solution to \eqref{0Int1} with $u\leq A$.
Denote by $\tilde{u}=u/A$, $f=\log\tilde{u}\leq0$ and
$w=|\nabla\log(1-f)|^2$. Then, it holds
\begin{equation}\label{Lemma1}
(\Delta-\partial_t) w\geq\frac{2f}{1-f}\nabla f\nabla w+2(1-f)w^2
-\frac{2(a+B)}{1-f}w,
\end{equation}
where $B=b+a\log A$.

\end{lem}

\proof Under the scaling $u\rightarrow \tilde{u}=u/A$, we have
$0<\tilde{u}\leq 1$. From \eqref{0Int1}, we obtain that $\tilde{u}$
satisfies the following equation
\begin{equation}\label{Proof1}
\tilde{u}_t=\Delta \tilde{u}+a\tilde{u}\log \tilde{u}+B\tilde{u},
\end{equation} where the constant $B$ satisfies $B=b+a\log A$. Let
$f=\log\tilde{u}\leq0$ and $w=|\nabla\log(1-f)|^2=\frac{|\nabla
f|^2}{(1-f)^2}$. Then we have
\begin{equation}\label{Proof2}
f_t=\Delta f+|\nabla f|^2+af+B.
\end{equation} Using \eqref{1Int3}, we can obtain
\begin{equation}\label{Proof4}
(|\nabla f|^2)_t=2R_{ij}f_{i}f_{j}+2f_i(f_t)_i.
\end{equation}
By the definition of $w$, we have
\begin{equation}\label{Proof5}\aligned
w_t=&\frac{2}{(1-f)^2}[R_{ij}f_{i}f_{j}+f_i(f_t)_i]
+\frac{2}{(1-f)^3}f_j^2f_t\\
=&\frac{2}{(1-f)^2}[R_{ij}f_{i}f_{j}+f_if_{jji}
+2f_{ij}f_if_j+af_i^2]\\
&+\frac{2}{(1-f)^3}f_j^2(f_{ii}+f_i^2+af+B).
\endaligned\end{equation}
On the other hand,
\begin{equation}\label{Proof6}\aligned
\Delta w=&\frac{6}{(1-f)^{4}}f_{i}^2f_{j}^2+\frac{2}{(1-f)^{3}}f_{ii}f_{j}^2
+\frac{8}{(1-f)^{3}}f_{ji}f_{i}f_{j}\\
&+\frac{2}{(1-f)^{2}}f_{ji}^2+\frac{2}{(1-f)^{2}}f_{j}f_{jii}\\
=&\frac{6}{(1-f)^{4}}f_{i}^2f_{j}^2+\frac{2}{(1-f)^{3}}f_{ii}f_{j}^2
+\frac{8}{(1-f)^{3}}f_{ji}f_{i}f_{j}\\
&+\frac{2}{(1-f)^{2}}f_{ji}^2+\frac{2}{(1-f)^{2}}f_{j}f_{iij}+\frac{2}{(1-f)^{2}}R_{ij}f_{i}f_{j},
\endaligned
\end{equation}
where, in the second equality, we used the Ricci
formula:$$f_{jii}=f_{iji}=f_{iij}+R_{ij}f_{i}.$$ By the formulas
\eqref{Proof5} and \eqref{Proof6}, we arrive at
\begin{equation}\label{Proof7}\aligned
(\Delta-\partial_t)
w=&\frac{2}{(1-f)^{2}}f_{ji}^2+\frac{8}{(1-f)^{3}}f_{ji}f_{i}f_{j}-\frac{4}{(1-f)^2}f_{ij}f_if_j
\\
&+\frac{6}{(1-f)^{4}}f_{i}^2f_{j}^2-\frac{2}{(1-f)^3}f_{i}^2f_{j}^2
-\frac{2(a+B)}{1-f}\frac{f_i^2}{(1-f)^2}.
\endaligned
\end{equation}
Note that
\begin{equation}\label{Proof8}
\nabla f\nabla
w=\frac{2}{(1-f)^{3}}f_{i}^2f_{j}^2+\frac{2}{(1-f)^{2}}f_{ji}
f_{i}f_{j}.
\end{equation} Therefore, \eqref{Proof7} can be written as
\begin{equation*}\label{Proof9}\aligned
(\Delta-\partial_t)
w=&\frac{2}{(1-f)^{2}}f_{ji}^2+\frac{4}{(1-f)^{3}}f_{ji}f_{i}f_{j}-\frac{4}{(1-f)^2}f_{ij}f_if_j
\\
&+\frac{2}{1-f}\nabla f\nabla w+2\frac{|\nabla
f|^4}{(1-f)^{4}}-2\frac{|\nabla f|^4}{(1-f)^3}
-\frac{2(a+B)}{1-f}\frac{|\nabla f|^2}{(1-f)^2}\\
=&\frac{2}{(1-f)^{2}}f_{ji}^2+\frac{4}{(1-f)^{3}}f_{ji}f_{i}f_{j}-2\nabla
f\nabla w
\\
&+\frac{2}{1-f}\nabla f\nabla w+2\frac{|\nabla
f|^4}{(1-f)^{4}}+2\frac{|\nabla f|^4}{(1-f)^3}
-\frac{2(a+B)}{1-f}\frac{|\nabla f|^2}{(1-f)^2}\\
=&\frac{2}{(1-f)^{2}}\left(f_{ji}+\frac{1}{1-f}f_{i}f_{j}\right)^2
\\
&+\frac{2f}{1-f}\nabla f\nabla w+2\frac{|\nabla f|^4}{(1-f)^3}
-\frac{2(a+B)}{1-f}\frac{|\nabla f|^2}{(1-f)^2}\\
\geq&\frac{2f}{1-f}\nabla f\nabla w+2(1-f)w^2 -\frac{2(a+B)}{1-f}w.
\endaligned
\end{equation*} Then, the desired estimate \eqref{Lemma1} follows.
\endproof

\noindent{\bf Proof of Theorem \ref{thmInt1}.} Let $\psi=\psi(x,t)$
be a smooth cutoff function supported in $B_{R,T}$ which satisfies
the following properties (see Page 13 in \cite{Zhang2006}, or Lemma
2.1 in \cite{Cao2010}):

1) $\psi=\psi(d(x,x_0,t),t)\equiv\psi(r,t)$; $\psi(x,t)=1$ in
$B_{\frac{R}{2},\frac{T}{2}}$ and $0\leq\psi\leq1$.

2) $\psi$ is decreasing as a radial function in the spatial
variables.

3) $|\frac{\partial\psi}{\partial
r}|\leq\frac{C_\gamma}{R}\psi^\gamma$ and
$|\frac{\partial^2\psi}{\partial
r^2}|\leq\frac{C_\gamma}{R^2}\psi^\gamma$ for every $\gamma\in
(0,1)$.

4) $|\frac{\partial\psi}{\partial t}|\leq\frac{C
}{T}\psi^{\frac{1}{2}}$.

We use \eqref{Lemma1} to conclude that
\begin{equation}\label{Proof13}\aligned
(\Delta-\partial_t) (\psi
w)=&\psi(\Delta-\partial_t)w+w(\Delta-\partial_t)\psi+2\nabla\psi\nabla
w\\
=&\psi(\Delta-\partial_t)w+w(\Delta-\partial_t)\psi\\
&+2\frac{\nabla\psi}{\psi}\nabla (\psi
w)-2w\frac{|\nabla\psi|^2}{\psi}\\
\geq&\frac{2f}{1-f}[\nabla f\nabla(\psi w)-w\nabla f\nabla
\psi]+2(1-f)\psi w^2\\
&-\frac{2(a+B)}{1-f}\psi
w+w(\Delta-\partial_t)\psi\\
&+2\frac{\nabla\psi}{\psi}\nabla (\psi
w)-2w\frac{|\nabla\psi|^2}{\psi}.
\endaligned\end{equation}
Now we let $(x_1,t_1)$ be a maximum point of $\psi w$ in the closure
of $B_{R,T}$, and where $\psi w>0$ (otherwise the proof is trivial).
Then at the point $(x_1,t_1)$, we have $\Delta(\psi w)\leq0$, $(\psi
w)_t\geq0$ and $\nabla(\psi w)=0$. Thus, from \eqref{Proof13}, we
deduce
\begin{equation}\label{Proof14}\aligned
2(1-f)\psi w^2\leq&\frac{2f}{1-f}w\nabla f\nabla \psi
+\frac{2(a+B)}{1-f}\psi
w\\
&-w(\Delta-\partial_t)\psi +2w\frac{|\nabla\psi|^2}{\psi}
\endaligned\end{equation} at $(x_1,t_1)$. Next, we will find an
upper bound for each term of the right hand side of \eqref{Proof14}.
By the Cauchy inequality, we have
\begin{equation}\label{Proof15}\aligned
\frac{2f}{1-f}w\nabla f\nabla \psi\leq&2|f|w^{\frac{3}{2}}\,|\nabla \psi|\\
=&2[\psi(1-f)w^2]^{\frac{3}{4}}\,\frac{f|\nabla \psi|}{[\psi(1-f)]^{\frac{3}{4}}}\\
\leq&(1-f)\psi w^2+C\frac{1}{R^4}\frac{f^4}{(1-f)^3}.
\endaligned\end{equation}
It has been shown in \cite{Zhang2006}(see formulas (3.30), (3.32)
and (3.34) in \cite{Zhang2006}) that:
\begin{equation}\label{Proof16}
w\frac{|\nabla\psi|^2}{\psi}\leq\frac{1}{8}\psi
w^2+C\frac{1}{R^4},\end{equation}
\begin{equation}\label{Proof17}\aligned
-(\Delta\psi)w\leq\frac{1}{8}\psi w^2+C\frac{1}{R^4}+C\frac{k}{R^2}
\endaligned\end{equation} and
\begin{equation}\label{Proof18}\aligned
\left|\frac{\partial\psi}{\partial t}\right|w\leq\frac{1}{8}\psi
w^2+C\frac{1}{T^2}+Ck^2.
\endaligned\end{equation}
Putting \eqref{Proof15}-\eqref{Proof18} into \eqref{Proof14}, we
obtain
\begin{equation}\label{Proof19}\aligned
(1-f)\psi w^2\leq&C\frac{1}{R^4}\frac{f^4}{(1-f)^3}
+\frac{2(a+B)}{1-f}\psi
w\\
&+\frac{1}{2}\psi w^2+C\frac{1}{R^4}+C\frac{1}{T^2}+Ck^2.
\endaligned\end{equation}
Since $f\leq0$, the inequality \eqref{Proof19} implies
\begin{equation}\label{Proof20}\aligned
\psi
w^2\leq&C\Bigg(\frac{1}{R^4}\frac{f^4}{(1-f)^4}+\frac{1}{R^4}+\frac{1}{T^2}+k^2
+(\max\{0,a+B\})^2\Bigg)\\
=&C\Bigg(\frac{1}{R^4}\frac{f^4}{(1-f)^4}+\frac{1}{R^4}+\frac{1}{T^2}+k^2+M_{a,b}^2\Bigg).
\endaligned\end{equation}
Therefore, we deduce that, for any $(x,t)\in B_{R,T}$,
\begin{equation}\label{Proof21}\aligned
(\psi^2 w^2)(x,t)\leq&(\psi
w^2)(x_1,t_1)\\
\leq&C\Bigg(\frac{1}{R^4}+\frac{1}{T^2}+k^2+M_{a,b}^2\Bigg)
\endaligned\end{equation} from $\frac{f^4}{(1-f)^4}<1$. Noticing
$\psi(x,t)=1$ in $B_{\frac{R}{2},\frac{T}{2}}$ and $w=\frac{|\nabla
f|^2}{(1-f)^2}$, we have
\begin{equation}\label{Proof22}
\frac{|\nabla f|}{1-f} \leq
C\Bigg(\frac{1}{R}+\frac{1}{\sqrt{T}}+\sqrt{k}+\sqrt{M_{a,b}}\Bigg)
\end{equation} holds at any $(x,t)\in
B_{\frac{R}{2},T}$, which shows that
\begin{equation}\label{Proof23}
\frac{|\nabla u|}{u} \leq
C\Bigg(\frac{1}{R}+\frac{1}{\sqrt{T}}+\sqrt{k}+\sqrt{M_{a,b}}\Bigg)\Bigg(1+\log\frac{A}{u}\Bigg).
\end{equation}
We complete the proof of the estimate \eqref{1thmFormula1} in
Theorem \ref{thmInt1}.

\section{Proof of Theorem \ref{thmInt2}}

As in the proof of Theorem \ref{thmInt1}, we first give a key lemma:

\begin{lem}\label{3Lem1}
Let $M$ be a complete Riemannian manifold with a family of
Riemannian metric $g(t)$ evolving by the Ricci flow \eqref{1Int3}.
Let $u$ be a positive solution to \eqref{1Int2} with $u\leq B$.
Denote by $\tilde{u}=u/B$, $f=\log\tilde{u}\leq0$ and
$w=|\nabla\log(1-f)|^2$. Then, it holds
\begin{equation}\label{3Lemma1}
(\Delta-\partial_t) w\geq\frac{2f}{1-f}\langle\nabla f,\nabla
w\rangle+2(1-f)w^2-2\tilde{\lambda}\Big(\alpha-\frac{-f}{1-f}\Big)e^{(\alpha-1)f}w,
\end{equation}
where $\tilde{\lambda}=\lambda B^{\alpha-1}$.

\end{lem}

\proof By the scaling $u\rightarrow \tilde{u}=u/B$, we have
$0<\tilde{u}\leq1$. Therefore, we obtain from \eqref{1Int2} that
$\tilde{u}$ satisfies
\begin{equation}\label{2Proof1}
\tilde{u}_t=\Delta \tilde{u}+\tilde{\lambda}\tilde{u}^{\alpha},
\end{equation} where $\tilde{\lambda}=\lambda B^{\alpha-1}$.
Let $f=\log \tilde{u}\leq0$ and
\begin{equation}\label{add2Proof1}w=|\nabla \log(1-f)|^2.\end{equation}
Then, the function $f$ satisfies
\begin{equation}\label{2Proof2}
f_t=\Delta f+|\nabla f|^2+\tilde{\lambda} e^{(\alpha-1)f}.
\end{equation}
Using \eqref{1Int3} again, one has
\begin{equation}\label{2Proof3}
(|\nabla f|^2)_t=2R_{ij}f_{i}f_{j}+2f_i(f_t)_i.
\end{equation}
It follows from \eqref{add2Proof1} that
\begin{equation}\label{2Proof4}\aligned
w_t=&\frac{2}{(1-f)^2}[R_{ij}f_{i}f_{j}+f_i(f_t)_i]
+\frac{2}{(1-f)^3}f_j^2f_t\\
=&\frac{2}{(1-f)^2}[R_{ij}f_{i}f_{j}+f_if_{jji}
+2f_{ij}f_if_j+\tilde{\lambda}(\alpha-1)e^{(\alpha-1)f}f_{i}^2]\\
&+\frac{2}{(1-f)^3}f_j^2(f_{ii}+f_i^2+\tilde{\lambda}
e^{(\alpha-1)f}).
\endaligned\end{equation}
Similarly, by the Ricci formula, we obtain
\begin{equation}\label{2Proof5}\aligned
\Delta
w=&\frac{6}{(1-f)^{4}}f_{i}^2f_{j}^2+\frac{2}{(1-f)^{3}}f_{ii}f_{j}^2
+\frac{8}{(1-f)^{3}}f_{ji}f_{i}f_{j}\\
&+\frac{2}{(1-f)^{2}}f_{ji}^2+\frac{2}{(1-f)^{2}}f_{j}f_{iij}+\frac{2}{(1-f)^{2}}R_{ij}f_{i}f_{j}.
\endaligned
\end{equation} Thus, we derive from \eqref{2Proof4} and  \eqref{2Proof5}
\begin{equation}\label{2Proof6}\aligned
(\Delta-\partial_t) w=&\frac{2}{(1-f)^{2}}f_{ji}^2+\Big[\frac{6}{(1-f)^{4}}-\frac{2}{(1-f)^{3}}\Big]f_{i}^2f_{j}^2\\
&+\Big[\frac{8}{(1-f)^{3}}-\frac{4}{(1-f)^{2}}\Big]f_{ji}f_{i}f_{j}\\
&-2\tilde{\lambda}\Big[(\alpha-1)\frac{1}{(1-f)^{2}}+\frac{1}{(1-f)^{3}}\Big]e^{(\alpha-1)f}f_{i}^2.
\endaligned
\end{equation}
Using the relationship
\begin{equation}\label{2Proof7}
\langle\nabla f,\nabla
w\rangle=\frac{2}{(1-f)^{3}}f_{i}^2f_{j}^2+\frac{2}{(1-f)^{2}}f_{ji}
f_{i}f_{j},
\end{equation} \eqref{2Proof6} can be written as
\begin{equation}\label{2Proof9}\aligned
(\Delta-\partial_t)& w-\varepsilon\langle\nabla f,\nabla
w\rangle\\
=&\frac{2}{(1-f)^{2}}f_{ji}^2+2\Big[\frac{3}{(1-f)^{4}}-(1+\varepsilon)\frac{1}{(1-f)^{3}}\Big]f_{i}^2f_{j}^2\\
&+2\Big[\frac{4}{(1-f)^{3}}-(2+\varepsilon)\frac{1}{(1-f)^{2}}\Big]f_{ji}f_{i}f_{j}\\
&-2\tilde{\lambda}\Big[(\alpha-1)\frac{1}{(1-f)^{2}}+\frac{1}{(1-f)^{3}}\Big]e^{(\alpha-1)f}f_{i}^2,
\endaligned
\end{equation} where $\varepsilon=\varepsilon(f)$ is a function depending on $f$ which will be
determined. Applying the inequality
$$\aligned
\frac{2}{(1-f)^{2}}&f_{ji}^2+2\Big[\frac{4}{(1-f)^{3}}-(2+\varepsilon)\frac{1}{(1-f)^{2}}\Big]f_{ji}f_{i}f_{j}\\
=&\frac{2}{(1-f)^{2}}\bigg\{f_{ji}^2+\Big[\frac{4}{1-f}-(2+\varepsilon)\Big]f_{ji}f_{i}f_{j}\bigg\}\\
\geq&-\frac{1}{2(1-f)^{2}}\Big[\frac{4}{1-f}-(2+\varepsilon)\Big]^2f_{i}^2f_{j}^2,
\endaligned$$ into \eqref{2Proof9} gives
\begin{equation}\label{2Proof10}\aligned
(\Delta-\partial_t)& w-\varepsilon\langle\nabla f,\nabla
w\rangle\\
\geq&\Big[-\frac{2}{(1-f)^{4}}+(6+2\varepsilon)\frac{1}{(1-f)^{3}}-(2+\varepsilon)^2\frac{1}{2(1-f)^{2}}\Big]f_{i}^2f_{j}^2\\
&-2\tilde{\lambda}\Big[(\alpha-1)\frac{1}{(1-f)^{2}}
+\frac{1}{(1-f)^{3}}\Big]e^{(\alpha-1)f}f_{i}^2\\
=&\frac{1}{(1-f)^{4}}\bigg\{-\frac{1}{2}(1-f)^{2}\varepsilon^2-2[(1-f)^{2}-(1-f)]\varepsilon\\
&-[2(1-f)^{2}-6(1-f)+2]\bigg\}f_{i}^2f_{j}^2\\
&-2\tilde{\lambda}\Big[(\alpha-1)\frac{1}{(1-f)^{2}}+\frac{1}{(1-f)^{3}}\Big]e^{(\alpha-1)f}f_{i}^2.
\endaligned
\end{equation}
Taking $$\varepsilon=-2+\frac{2}{1-f}$$ in \eqref{2Proof10}, we
derive
\begin{equation}\label{2Proof11}\aligned
(\Delta-\partial_t)& w-\frac{2f}{1-f}\langle\nabla f,\nabla
w\rangle\\
\geq&\frac{2}{(1-f)^{3}}|\nabla
f|^4-2\tilde{\lambda}\Big[(\alpha-1)\frac{1}{(1-f)^{2}}+\frac{1}{(1-f)^{3}}\Big]e^{(\alpha-1)f}f_{i}^2\\
=&2(1-f)w^2-2\tilde{\lambda}\Big(\alpha-\frac{-f}{1-f}\Big)e^{(\alpha-1)f}w.
\endaligned
\end{equation} Thus, the desired estimate \eqref{3Lemma1} is
attained.\endproof

\noindent{\bf Proof of Theorem \ref{thmInt2}.} Let $\psi=\psi(x,t)$
be a smooth cutoff function supported in $B_{R,T}$ which is defined
in section 2. We use \eqref{3Lemma1} to conclude that
\begin{equation}\label{2Proof12}\aligned
(\Delta-\partial_t) (\psi
w)=&\psi(\Delta-\partial_t)w+w(\Delta-\partial_t)\psi\\
&+2\frac{\nabla\psi}{\psi}\nabla (\psi
w)-2w\frac{|\nabla\psi|^2}{\psi}\\
\geq&-2\frac{-f}{1-f}[\nabla f\nabla(\psi w)-w\nabla f\nabla
\psi]+2(1-f)\psi w^2\\
&-2\tilde{\lambda}\Big(\alpha-\frac{-f}{1-f}\Big)e^{(\alpha-1)f}\psi
w+w(\Delta-\partial_t)\psi\\
&+2\frac{\nabla\psi}{\psi}\nabla (\psi
w)-2w\frac{|\nabla\psi|^2}{\psi}.
\endaligned\end{equation}
Now we let $(x_2,t_2)$ be a maximum point of $\psi w$ in the closure
of $B_{R,T}$, and where $\psi w>0$ (otherwise the proof is trivial).
Then at the point $(x_2,t_2)$, we have $\Delta(\psi w)\leq0$, $(\psi
w)_t\geq0$ and $\nabla(\psi w)=0$. Thus, from \eqref{2Proof12}, we
deduce
\begin{equation}\label{2Proof13}\aligned
2(1-f)\psi w^2\leq&\frac{2f}{1-f}w\nabla f\nabla \psi
+2\tilde{\lambda}\Big(\alpha-\frac{-f}{1-f}\Big)e^{(\alpha-1)f}\psi
w\\
&-w(\Delta-\partial_t)\psi +2w\frac{|\nabla\psi|^2}{\psi}
\endaligned\end{equation} at $(x_2,t_2)$.

{\bf Case one:} $\alpha\geq1$.

Note that
\begin{equation}\label{2Proof14}\frac{-f}{1-f}=1-\frac{1}{1-f}\in(0,1).\end{equation}
In this case, we have $\alpha-\frac{-f}{1-f}>0$ and
$e^{(\alpha-1)f}\in(0,1)$. Hence
\begin{equation}\label{2Proof15}
2\tilde{\lambda}\Big(\alpha-\frac{-f}{1-f}\Big)e^{(\alpha-1)f}\psi
w\leq
C\frac{M_{\lambda}^2}{1-f}\Big(\alpha-\frac{-f}{1-f}\Big)^2+\frac{1}{16}(1-f)\psi
w^2,\end{equation} where $M_{\lambda}=\max\{0,\tilde{\lambda}\}$.
Putting \eqref{2Proof15}, \eqref{Proof15}-\eqref{Proof18} into
\eqref{2Proof13}, we obtain
\begin{equation}\label{2Proof16}\aligned
(1-f)\psi w^2\leq&C\Bigg[\frac{1}{R^4}\frac{f^4}{(1-f)^3}
+\frac{M_{\lambda}^2}{1-f}\Big(\alpha-\frac{-f}{1-f}\Big)^2\\
&+\frac{1}{R^4}+\frac{1}{T^2}+k^2\Bigg].
\endaligned\end{equation}
Therefore, we deduce, for any $(x,t)\in B_{R,T}$,
\begin{equation}\label{2Proof18}\aligned
(\psi^2 w^2)(x,t)\leq&(\psi
w^2)(x_1,t_1)\\
\leq&C\Bigg(\frac{1}{R^4}+\frac{1}{T^2}+k^2+M_{\lambda}^2\alpha^2\Bigg)
\endaligned\end{equation} from $\frac{f^2}{(1-f)^2}<1$. Noticing
$\psi(x,t)=1$ in $B_{\frac{R}{2},\frac{T}{2}}$ and $w=\frac{|\nabla
f|^2}{(1-f)^2}$, we have
\begin{equation}\label{2Proof19}
\frac{|\nabla f|}{1-f} \leq
C\Bigg(\frac{1}{R}+\frac{1}{\sqrt{T}}+\sqrt{k}+\sqrt{M_{\lambda}\alpha}\Bigg)
\end{equation} holds at any $(x,t)\in
B_{\frac{R}{2},T}$, which shows that
\begin{equation}\label{2Proof20}
\frac{|\nabla u|}{u} \leq
C\Bigg(\frac{1}{R}+\frac{1}{\sqrt{T}}+\sqrt{k}
+\sqrt{M_{\lambda}\alpha}\Bigg)\Bigg(1+\log\frac{B}{u}\Bigg).
\end{equation}

{\bf Case two:} $\alpha\leq0$.

In this case, we have $\alpha-\frac{-f}{1-f}<0$ and
$e^{(\alpha-1)f}>1$. Hence
\begin{equation}\label{2Proof21}
2\tilde{\lambda}\Big(\alpha-\frac{-f}{1-f}\Big)e^{(\alpha-1)f}\psi
w\leq
C\frac{\tilde{M}_{\lambda}^2}{1-f}\Big(\alpha-\frac{-f}{1-f}\Big)^2\tilde{u}_{\min}^{2(\alpha-1)}
+\frac{1}{16}(1-f)\psi w^2,\end{equation} where
$\tilde{M}_{\lambda}=\max\{0,-\tilde{\lambda}\}$. Similarly, putting
\eqref{2Proof21}, \eqref{Proof15}-\eqref{Proof18} into
\eqref{2Proof13}, we obtain
\begin{equation}\label{2Proof22}\aligned
(1-f)\psi w^2\leq&C\Bigg[\frac{1}{R^4}\frac{f^4}{(1-f)^3}
+\frac{\tilde{M}_{\lambda}^2}{1-f}\Big(\alpha-\frac{-f}{1-f}\Big)^2\tilde{u}_{\min}^{2(\alpha-1)}\\
&+\frac{1}{R^4}+\frac{1}{T^2}+k^2\Bigg].
\endaligned\end{equation} Thus, we can obtain
\begin{equation}\label{2Proof23}
\frac{|\nabla f|}{1-f} \leq
C\Bigg(\frac{1}{R}+\frac{1}{\sqrt{T}}+\sqrt{k}
+\sqrt{\tilde{M}_{\lambda}(-\alpha+1)\tilde{u}_{\min}^{\alpha-1}}\Bigg)
\end{equation} holds at any $(x,t)\in
B_{\frac{R}{2},T}$, which shows that
\begin{equation}\label{2Proof24}
\frac{|\nabla u|}{u} \leq
C\Bigg(\frac{1}{R}+\frac{1}{\sqrt{T}}+\sqrt{k}
+\sqrt{M_{\lambda}(-\alpha+1)u_{\min}^{\alpha-1}}\Bigg)\Bigg(1+\log\frac{B}{u}\Bigg).
\end{equation}

{\bf Case three:} $\alpha\in(0,1)$.

In this case, we also have $e^{(\alpha-1)f}>1$. Hence
\begin{equation}\label{2Proof25}
2\tilde{\lambda}\Big(\alpha-\frac{-f}{1-f}\Big)e^{(\alpha-1)f}\psi
w\leq
C\frac{\tilde{\lambda}^2}{1-f}\Big(\alpha-\frac{-f}{1-f}\Big)^2\tilde{u}_{\min}^{2(\alpha-1)}
+\frac{1}{16}(1-f)\psi w^2.\end{equation} Putting \eqref{2Proof25},
\eqref{Proof15}-\eqref{Proof18} into \eqref{2Proof13}, we obtain
\begin{equation}\label{2Proof26}\aligned
(1-f)\psi w^2\leq&C\Bigg[\frac{1}{R^4}\frac{f^4}{(1-f)^3}
+\frac{\tilde{\lambda}^2}{1-f}\Big(\alpha-\frac{-f}{1-f}\Big)^2\tilde{u}_{\min}^{2(\alpha-1)}\\
&+\frac{1}{R^4}+\frac{1}{T^2}+k^2\Bigg].
\endaligned\end{equation} Because of $\Big(\alpha-\frac{-f}{1-f}\Big)^2\leq1$, we can obtain
\begin{equation}\label{2Proof27}
\frac{|\nabla f|}{1-f} \leq
C\Bigg(\frac{1}{R}+\frac{1}{\sqrt{T}}+\sqrt{k}
+\sqrt{|\tilde{\lambda}|\tilde{u}_{\min}^{\alpha-1}}\Bigg)
\end{equation} holds at any $(x,t)\in
B_{\frac{R}{2},T}$, which shows that
\begin{equation}\label{2Proof28}
\frac{|\nabla u|}{u} \leq
C\Bigg(\frac{1}{R}+\frac{1}{\sqrt{T}}+\sqrt{k}
+\sqrt{|\lambda|u_{\min}^{\alpha-1}}\Bigg)\Bigg(1+\log\frac{B}{u}\Bigg).
\end{equation}
We complete the proof of Theorem \ref{thmInt2}.
%
%

\bibliographystyle{Plain}

\end{document}